\documentclass[12pt,reqno]{amsart}
\usepackage{amssymb,array}

\newcommand{\PBS}[1]{\let\temp=\\#1\let\\=\temp}

\newtheorem{lemma}{Lemma}
\newtheorem{proposition}[lemma]{Proposition}
\newtheorem{theorem}[lemma]{Theorem}
\newtheorem{corollary}[lemma]{Corollary}

\newtheorem{comm}[lemma]{Comment}
\newtheorem{remm}[lemma]{Remark}
\newtheorem{exam}[lemma]{Example}

\newenvironment{vect}{\left( \!\!\! \begin{array}{c}}{\end{array} \!\!\!\right)}
\newenvironment{matr2}{\left( \!\!\! \begin{array}{cc}}{\end{array} \!\!\!\right)}

\newcommand{\alp}{\alpha}
\newcommand{\bet}{\beta}

\newcommand{\lam}{\lambda}

\newcommand{\del}{\delta}
\newcommand{\eps}{\varepsilon}
\newcommand{\tet}{\vartheta}

\newcommand{\Lam}{\Lambda}

\newcommand{\LL}{\mathcal{L}}

\newcommand{\R}{\mathbb{R}}
\newcommand{\C}{\mathbb{C}}

\newcommand{\Z}{\mathbb{Z}}
\newcommand{\T}{\mathbb{T}}
\newcommand{\Le}[1][2]{{\rm L}^{#1}}
\renewcommand{\le}[1][2]{{\rm l}^{#1}}

\newcommand{\ud}{{\rm d}}
\newcommand{\e}{{\rm e}}

\newcommand{\hull}{\mathrm{\hull}}

\newcommand{\spec}[2][\!]{\mathrm{Spec}_{#1} \left( #2 \right) \, }

\newcommand{\dist}[1]{\mathrm{dist} \left( #1 \right) \, }

\newcommand{\dotp}[1]{\left< #1 \right>}

\newcommand{\trac}[1]{\mathrm{Tr} \left( #1 \right) \, }

\newcommand{\Proof}{\textbf{\underline{Proof}.}\ }
\newcommand{\fin}{\hspace*{\fill} $\, \blacksquare$}
\newcommand{\find}{\hspace{.1in} \blacksquare}

\fontshape{n}

\title[Projection methods...]
{Projection methods for discrete Schr\"odinger operators}
\author[L. Boulton]{}

\date{January 2003}



\thanks{$^1$ The author carried out part of this research at
King's College London, founded by EPSRC under grant number
GR/L75443.}
\begin{document}
\

%%%%%%%%%%%%%%%%%%%%%%%%%%%%%%%%%%%%%%%
% PAPER NUMBER
%\vskip -1in {\em \
%\parbox[t]{1.1\linewidth}{\small
%\rightline{number} } }

%\vskip -1.5em \null

%\vskip 1.5em
%%%%%%%%%%%%%%%%%%%%%%%%%%%%%%%%%%%%%%%%

\maketitle

\centerline{\scshape L.~Boulton$^1$}

\medskip

{\footnotesize \centerline{Departmento de Matem\'aticas,
Universidad Sim\'on Bol\'\i var} \centerline{Apartado 89000,
Caracas 1080-A, Venezuela} \centerline{email: lboulton@ma.usb.ve}}

\medskip

\begin{abstract}
Let $H$ be the discrete Schr\"odinger operator \linebreak
$Hu(n):=u(n-1)+u(n+1)+v(n)u(n)$, $u(0)=0$ acting on $\le (\Z^+)$
where the potential $v$ is real-valued and $v(n)\to 0$ as $n\to
\infty$. Let $P$ be the orthogonal projection onto a closed linear
subspace $\LL \subset \le (\Z^+)$. In a recent paper E.B. Davies
defines the second order spectrum ${\rm Spec}_2(H,\LL)$ of $H$
relative to $\LL$ as the set of $z \in \C$ such that the
restriction to $\LL$ of the operator $P(H-z)^2P$ is not invertible
within the space $\LL$. The purpose of this article is to
investigate properties of ${\rm Spec}_2(H,\LL)$ when $\LL$ is
large but finite dimensional. We explore in particular the
connection between this set and the spectrum of $H$. Our main
result provides sharp bounds in terms of the potential $v$ for the
asymptotic behaviour of ${\rm Spec}_2(H,\LL)$ as $\LL$ increases
towards $\le(\Z^+)$.
\end{abstract}

\section{Introduction} \label{s1}
Let the discrete Schr\"odinger operator $H$ be defined by
\[
   Hu(n):=u(n+1)+u(n-1)+v(n)u(n),\qquad u(0)=0
\]
acting on $\le (\Z^+)$ where $v:\Z^+ \longrightarrow \R$.
In \cite{SECR}, E.B.~Davies investigated the concept of resonance
of $H$ associated to a closed subspace \linebreak  $\LL\subset \le (\Z^+)$.
Let $P$ be the orthogonal projection onto $\LL$.
He called the isolated points of the second order spectrum of $H$
\begin{equation} \label{e19}
  \spec[2]{H;\LL} := \{ \lam \in \C \,:\, P(H-\lam)^2P|\LL\
   \mathrm{is\ not\ invertible}\},
\end{equation}
resonances of $H$ relative to $\LL$. This definition is motivated
by graphical and numerical connections between other notions of
resonance adopted to this context, and isolated points of
$\spec[2]{H;\LL}$ (cf. \cite[section 9]{SECR}). This was observed
by Davies in the case when the dimension of $\LL$ is low ($\leq
30$) and $v$ has finite support.

In \cite{SECR} Davies also pointed out that $\spec[2]{H;\LL}$
contains information about the spectrum $\spec{H}$ of $H$.
 In \cite{GHORS}
E.~Shargorodsky developed this idea for general linear operators
$T$. This began the study of second order spectra as a projection
method for localizing the spectrum. Among other interesting
results on geometrical properties of $\mathrm{Spec}_2$,
Shargorodsky showed that for any bounded self-adjoint operator
$T$,
\begin{equation} \label{e18}
 \bigcup \lim_{k\to \infty} \spec[2]{T;\LL_k} \cap \R =\spec{T}
\end{equation}
where the union is taken over the set of all sequences $(\LL_k)$
of subspaces of the domain of $T$ such that the orthogonal
projection $P_k$ onto $\LL_k$ converges strongly to the identity,
see \cite[theorem 21]{GHORS}. Hence the second order spectra of
$H$ might be useful for approximate computation of $\spec{H}$ when
$v$ is bounded.

The purpose of the present paper is to continue investigating the
connection between the spectrum and second order spectra of $H$.
We present a detailed description of
\[
 \spec[2]{H;k}:=\spec[2]{H;\LL_k}
\]
for large $k$ when
\begin{equation*}
   \LL_k:=\mathrm{span} \{\del_1,\ldots,\del_k\}, \qquad
   \del_j(n):=\left\{\begin{array}{cc} 0 & n\not=j \\ 1 &n=j
 \end{array} \right.
\end{equation*}
and there exist constants $a,r>0$ such that
\begin{equation} \label{e2}
 |v(n)|\leq \frac{a}{n^{r+1}} \qquad \qquad n=1,2,\ldots
\end{equation}
As  $v$ is real-valued and bounded, $H$ is a bounded self-adjoint
operator so we are in the situation covered by Shargorodsky's
results. Since \eqref{e2} holds for finite rank $v$, our
discussion includes the potentials investigated by Davies in
\cite{SECR}.

By construction, $\LL_k$ is $k$-dimensional. Then
the quadratic operator pencil in \eqref{e19} is a $k\times k$ matrix and
\begin{equation*}
   \spec[2]{H;k} =  \{ \lam \in \C \,:\, \det[P_k(H-\lam)^2P_k|\LL_k]=0 \}.
\end{equation*}
It is easy to see that the determinant at the right hand side is
a polynomial in $\lam$ of degree at most $2k$ so that
$\spec[2]{H;k}$ consists exclusively of at most
$2k$ relative resonances of $H$.

Since our interest is in how the $\spec[2]{H;k}$ are related to
the spectrum of $H$, we describe the latter. The classical theory
allows us to find $\spec{H}$ under the hypothesis \eqref{e2} quite
efficiently. For this we decompose
\[
 H=H_0+V
\]
where $H_0$ is $H$ without potential and
$ Vu(n):= v(n)u(n)$. The spectrum of $H_0$ is pure
absolutely continuous and equal to the interval $[-2,2]$. Since
\begin{align*}
 \trac{V} & =
 \sum_{n=1}^\infty |\dotp{V\del_n,\del_n}|  = \sum_{n=1}^\infty
 |v(n)| \\ & \leq \sum_{n=1}^\infty \frac{a}{n^{r+1}} < +\infty,
\end{align*}
$V$ is a trace class operator. Then, by virtue of Pearson's theorem,

\[
  \spec[\mathrm{ess}]{H}=\spec[\mathrm{ac}]{H}=[-2,2]
\]
so $\spec{H}$ consists of the interval $[-2,2]$ together with a
discrete set of isolated eigenvalues of finite multiplicity which
can only accumulate at $\pm 2$. A recent adaptation by F.~Luef and
G.~Teschl  of the Sturm-Liouville theorem
to the discrete context, cf. \cite{LuTe}, shows that in fact
if $r>1$ in $\eqref{e2}$, then $H$ has
only a finite number of eigenvalues.

\medskip

The results we will discuss in the forthcoming sections are motivated by
the following proposition.

\begin{proposition} \label{t3}
Let $V=V^{\ast}$ be a compact operator acting on $\le (\Z^+)$.
Then for all $z\not\in \spec{H_0+V}$ there exists $\tilde{k}>0$
such that
\[
  z \not \in \spec[2]{H_0+V;k}, \qquad k\geq \tilde{k}.
\]
\end{proposition}

In other words, roughly speaking
\begin{equation} \label{e12}
 S:=\lim_{k\to \infty} \spec[2]{H;k} \subseteq \spec{H}.
\end{equation}
Notice that we are less restrictive on the requirements for $V$.
In comparing with \eqref{e18}, this result says that our
particular choice of $\LL_k$ does not lead to spurious points in
approximating $\spec{H}$. In section \ref{s2} we will show
proposition \ref{t3} using an argument that will be crucial for
our latter work. It consists in considering $\spec[2]{H;k}$ as the
perturbed spectrum of a convenient Toeplitz matrix.

Our main results are theorem~\ref{t7} and corollary~\ref{t8} of section~\ref{s4}.
There we give a precise meaning to the limit in \eqref{e12} and provide
sharp estimates on the rate of convergence in terms of $v$ when it satisfies
\eqref{e2}.

The proof of theorem~\ref{t7} requires to develop some results on
the stability of Toeplitz operators with smooth symbols. This is
carried out in section~\ref{s3}. It will become clear from the
beginning that these techniques can be extended to more general
situations, hence we give a general treatment to the topic.
Section~\ref{s3} is closely related to chapter~2 of the monograph
\cite{BoSi} on large Toeplitz matrices. Since we failed to find a
reference for theorem~\ref{t4} and corollary~\ref{t6} we include
proofs of both.

We complete our discussion by considering in detail some numerical examples in
section~\ref{s5}. These are closely related to those of \cite[section~9]{SECR}.

\section{The finite section method} \label{s2}

The finite section method for perturbed Toeplitz operators was
developed by  I.~Gohberg, I.A.~Feldman  and H.~Widom in the 1970's
(see  \cite{GF} and \cite{W}). Through this paper we follow
closely their ideas, in particular we follow the notation and
procedures of the introductory text on the subject by A.~B\"otcher
and B.~Silbermann, \cite{BoSi}. This section is devoted to discuss
the basic notation and results in the Gohberg-Feldman-Widom
approach which will be the core of our latter work. This will lead
us to a proof of proposition \ref{t3}.

Below and elsewhere we identify the truncation
\[
 P_kTP_k|\LL_k:
\LL_k \longrightarrow \LL_k
\]
of any linear operator $T:\le(\Z^+)\longrightarrow \le(\Z^+)$,
with its matrix representation acting on $\C^k$. We will often
denote $T_k:=P_kTP_k|\LL_k$. The following definition is taken
from \cite[chapter 2]{BoSi}: we say that the sequence
$\{P_kTP_k\}$ is stable (with constants $k_0$ and $M$) if and only
if, there exists $k_0>0$ such that $T_k$ is invertible for all $k
\geq k_0$ and
\begin{equation*}
 \sup_{k \geq k_0} \| T_k^{-1} \| \leq M
<\infty.
\end{equation*}
Notice that if $\{P_k(H-z)^2P_k\}$ is stable then $z\in \C$ falls
outside $\spec[2]{H;k}$ for large $k$.

If $T$ is a Toeplitz operator with continuous symbol, then
$\{P_kTP_k\}$ is stable if and only if $T$ is invertible (cf.
\cite[theorem 2.11]{BoSi}). Furthermore (see corollary~\ref{t6}
below or \cite[theorem 2.16]{BoSi}), if $\{P_kTP_k\}$ is stable,
$K$ is compact and $T+K$ is invertible, then $\{P_k(T+K)P_k\}$ is
stable.

\medskip

The following argument is crucial to out later work.  Decompose
\begin{equation} \label{e1} \begin{aligned}
 (H-z)^2 & = (H_0-z)^2+(H_0-z)V+V(H_0-z) +V^2 \\
 & = T(\bet)+K(z)
\end{aligned} \end{equation}
where $T(\bet)$ is the Toeplitz operator whose symbol is
\begin{align*}
 \bet(t;z)&:= t^{-2}-2zt^{-1}+(2+z^2)-2zt+t^2 \\
 & = (t^{-2}-zt^{-1}+1)(t^{2}-zt+1) \\
 & =(t^{-1}-z+t)^2,\\
 t\in & \,\T:=\{\e^{i\tet}:-\pi<\tet\leq \pi\}
\end{align*}
and
\begin{align*}
 K(z)& :=(H_0-z)V+V(H_0-z) +V^2-|\del_1\rangle \langle \del_1| \\
 & = H_0V+VH_0+V^2-2zV-|\del_1\rangle \langle \del_1|
\end{align*}
is compact. According to the above, the stability of
$\{P_kT(\bet)P_k\}$ implies stability of  $\{P_k(H-z)^2P_k\}$ when
$z\not \in \spec{H}$. This is the key argument in the proof of
proposition \ref{t3}.

\textbf{\underline{Proof} [of proposition \ref{t3}].} The symbol
\begin{equation*}
 \bet(\e^{i\tet};z)  = (2\cos\tet -z)^2, \qquad -\pi < \tet \leq \pi.
\end{equation*}
Then $\bet(\T;z)$ is an open curve without loops so therefore, by
virtue of Gohberg's theorem on the spectrum of Toeplitz operators
with continuous symbol (see for instance
\cite[theorem~1.17]{BoSi}), $\mathrm{Spec}\,T(\bet)=\bet(\T;z)$.
Thus $T(\bet)$ is invertible if and only if $z\not\in
[-2,2]$. For all such $z$, $\{P_kT(\bet)P_k\}$
is stable. Clearly $(H-z)^2$ is invertible if and only if
$z\not\in \spec{H}$. Hence,
\[
 \{P_k[T(\bet)+K(z)]P_k\}=\{P_k(H-z)^2P_k\}
\]
is stable for all $z \not \in \spec{H}$. \fin

\section{Stability of Toeplitz operators with smooth symbol} \label{s3}
In this section $T(\alp)$ denotes the Toeplitz operator associated
to the symbol $\alp:\T \longrightarrow \C$. If $q>0$, $W^q$ stands for the set of
continuous $\alp$ such that
\[
 \sum_{n\in\Z} |n|^q|\widehat{\alp}(n) |<\infty
\]
where $\widehat{\alp}(n)$ denote here and elsewhere the
(non-normalized) Fourier coefficients of $\alp$. It is well known
that $W^q$ is a Banach algebra with pointwise algebraic operations
and norm $\|\alp\|:=\sum (1+|n|)^q|\widehat{\alp}(n)|$. According
to Wiener's theorem (cf. \cite{wie}), if $\alp\in W^q$ and $0\not
\in \alp(\T)$, then \linebreak $\alp^{-1}\in W^q$.

Throughout this section we will assume that
\begin{itemize}
\item[i)] The symbol $\alp\in W^q$ for some $q>1$.
\item[ii)] The sequence $\{P_kT(\alp)P_k\}$ is stable with
constants denoted by $k_0$ and $M$: that is $T_k(\alp)$ is
invertible for all $k \geq k_0$ and \linebreak $\sup_{k \geq k_0}
\| T_k^{-1}(\alp) \| \leq M$.
\end{itemize}
The  continuity of $\alp$ and stability of $\{P_kT(\alp)P_k\}$
guarantee that $T(\alp)$ is invertible. Under these conditions it
is well known (cf. \cite[pp.32-34]{BoSi}) that $T^{-1}_k(\alp)$
converges in the strong operator topology to $T^{-1}(\alp)$.
Therefore it is legitimate to estimate the columns of the matrix
of $T^{-1}(\alp)$ from the columns of $T_k^{-1}(\alp)$ for large
$k$. The following result states that if we not only assume
continuity of the symbol but the stronger condition i), the error
in the estimation of the first $k/2$ columns of $T^{-1}(\alp)$
from $T_k^{-1}(\alp)$ is $O(k^{q-\eps})$.

\begin{theorem} \label{t4}
Let $T(\alp)$ satisfy i)-ii). Then for all $1\leq p< q$ there
exists a constant $c>0$ independent of $k$, such that
\[
 \|T_k^{-1}(\alp)\del_j-T^{-1}(\alp)\del_j\| \leq \frac{c}{k^p}
\]
provided $k\geq k_0$ and $1\leq j\leq k/2$.
\end{theorem}

Below we will compute an explicit upper bound for $c$ in terms of
$\alp$, $M$ and $\|T^{-1}(\alp)\|=:\tilde{M}$.

\medskip

This theorem improves \cite[theorem~2.15]{BoSi}. It allows us to
show a version of the result mention in section~\ref{s2} on
stability of compact perturbation of stable sequences, when the
perturbed operator $T$ is Toeplitz, its symbol satisfies i) and
the perturbation $K$ is a trace class band matrix. This is the
content of corollary \ref{t6}. Let
\begin{align*}
 K_{j,l}&=\dotp{K\del_l,\del_j}\qquad \qquad & l,j=1,2,\ldots \\
 K_{j,l}&=0 & l\leq0 \mathrm{\ or\ } j\leq0 .
\end{align*}
We assume that $K$ satisfies the following hypotheses,
\begin{itemize}
\item[a)] There exists $L> 0$ such that $K_{j,l}=0$
for $| j-l|> L$.
\item[b)] There exist $b>0$ and $r>0$, such that
\begin{equation*}
 K_{j,l}\leq \frac{b}{j^{r+1}},  \qquad \qquad l,j=1,2,\ldots
\end{equation*}
\end{itemize}

\begin{corollary} \label{t6}
Let $T(\alp)$ satisfy i)-ii) and let $K$ satisfy a)-b). Fix
\linebreak $1\leq p <q$, let $\rho:=\min\{r,(2p-1)/2\}$ and let
$c$ be as in theorem~\ref{t4}. If $T(\alp)+K$ is invertible, then
$P_k(T(\alp)+K)P_k|\LL_k$ is invertible provided
 \begin{align}
  k \geq & \nonumber \\ \, \max& \left\{
  \left(\frac{4L\left[c\sqrt{3}\|K\|+b 2^{r}\sqrt{\pi}(M+\tilde{M})\right]
  \|(T(\alp)+K)^{-1}T(\alp)\|}{\sqrt{6}}\right)^{1/\rho},k_0\right\}\label{e11}.
 \end{align}
\end{corollary}
\Proof Throughout the proof we assume $1\leq p<q$ and $k> k_0$. By
virtue of ii),
\[
 I+T^{-1}(\alp)K=T^{-1}(\alp)(T(\alp)+K)
\]
is invertible. Let
\[
 1/\eps:=\|(I+T^{-1}(\alp)K)^{-1}\|=\|(T(\alp)+K)^{-1}T(\alp)\|.
\]
Then for all $u\in \le (\Z^+)$,
\begin{align*}
 \eps \|P_ku \|&\leq \|(I+T^{-1}(\alp)K)P_ku\| \\
 & \leq \|(I+T_k^{-1}(\alp)P_kK)P_ku\| +
 \|T^{-1}_k(\alp)P_kK-T^{-1}(\alp)K\|
 \,\|P_ku\|.
\end{align*}
The operator $K$ acts on $u$ as follows,
\[
 Ku=\sum_{n=1}^\infty \sum _{m=-L}^{L} K_{n,n+m} u(n+m) \del_n.
\]
Let $R_k:=T^{-1}_k(\alp)P_k-T^{-1}(\alp)$. Then
\begin{align*}
 \|R_kKu\| & = \left\|\sum_{n=1}^\infty \sum _{m=-L}^{L} K_{n,n+m}
 u(n+m)R_k\del_n \right\| \\
 & \leq \sum_{m=-L}^{L} \sum_{n=1}^\infty
 |u(n+m)|\,|K_{n,n+m}|\,\|R_k \del_n\|.
\end{align*}
For each $m=-L,\ldots,L$,
\begin{align*}
 \sum_{n=1}^\infty &|u(n+m)|  \,|K_{n,n+m}|\,\|R_k \del_n\| \leq \\
 & \leq \left( \sum_{n=1}^\infty |u(n+m)|^2 \right)^{1/2}
 \left(\sum_{n=1}^\infty  |K_{n,n+m}|^2\|R_k \del_n\|^2\right)^{1/2}
 \\ & \leq \|u\| \left( \sum_{n=1}^{k/2}  |K_{n,n+m}|^2\|R_k \del_n\|^2
 +\sum_{n=k/2+1}^{\infty}  |K_{n,n+m}|^2\|R_k \del_n\|^2
 \right)^{1/2}.
\end{align*}
By virtue of theorem \ref{t4},
\begin{align*}
 \sum_{n=1}^{k/2}  |K_{n,n+m}|^2\|R_k \del_n\|^2  & \leq
 \frac{c^2\sum_{n=1}^{k/2}  |K_{n,n+m}|^2}{k^{2p}} \\
 & \leq \frac{c^2\|K\|^2\sum_{n=1}^{k/2}1}{k^{2p}} =
 \frac{c^2\|K\|^2}{2k^{2p-1}}.
\end{align*}
By virtue of a)-b) and by the definition of $R_k$,
\begin{align*}
 \sum_{n=k/2+1}^{\infty}  |K_{n,n+m}|^2\|R_k \del_n\|^2 & \leq
 \|R_k\|^2 \sum_{n=k/2+1}^{\infty} |K_{n,n+m}|^2 \\
 & \leq (M+\tilde{M})^2 \sum_{n=k/2+1}^{\infty} |K_{n,n+m}|^2 \\ &
 \leq \frac{b^22^{2r}(M+\tilde{M})^2 \sum_{n=k/2+1}^{\infty}
 \left[(k/2)^r/n^{r+1}\right]^2}{k^{2r}} \\ &
 \leq \frac{b^22^{2r}(M+\tilde{M})^2 \sum_{n=k/2+1}^{\infty}
 1/n^2}{k^{2r}} \\ & \leq
 \frac{b^22^{2r}(M+\tilde{M})^2 \pi
 }{6k^{2r}}.
\end{align*}
Then, if $k$ satisfies \eqref{e11},
\[
 \|T^{-1}_k(\alp)K-T^{-1}(\alp)K\|
  \leq  2L\frac{c\sqrt{3}\|K\|+b2^{r}\sqrt{\pi}
  (M+\tilde{M})}{\sqrt{6}k^{\rho}}\leq \eps /2.
\]
For all such $k$,
\begin{align*}
 (\eps/2)\|P_kv\| & \leq \|(I+T_k^{-1}(\alp)P_kK)P_kv\| \\
 & \leq \|T_k^{-1}(\alp)\|\,\|(T_k(\alp)+P_kK)P_kv\| \\
 & =\|T_k^{-1}(\alp)\|\,\|P_k(T(\alp)+K)P_kv\| \\
 & \leq M \|P_k(T(\alp)+K)P_kv\|,
\end{align*}
so that $P_k(T(\alp)+K)P_k|\LL_k$ is injective and thus
invertible. Notice that $\dim \LL_k$ is finite. \fin

In section \ref{s4} we will employ this corollary to estimate sharp
bounds for the error in the limit \eqref{e12}.

\medskip

We now prove theorem~\ref{t4}. For this purpose we use the
following convention: $Q_k:=I-P_k$,
\[
 W_ku(n):=\left\{
 \begin{array}{cc} u(k-n+1) & 1\leq n \leq k \\
 0 &n>k \end{array} \right.,
\]
$\tilde{\alp}:=\alp(t^{-1})$, $H(\alp)$ is the Hankel operator
generated by $\alp$ (cf. \cite[p.13]{BoSi}),
\[
 \Lam(\alp):=T^{-1}(\alp)-T(\alp^{-1})
\]
and
\[
  B_k:=P_kT^{-1}(\alp)P_k+W_k\Lam(\tilde{\alp})W_k.
\]
Assumption ii) and the fact that $\alp$ is continuous, ensure that
$\Lam(\alp)$ is well defined. The following results are due to
Widom, but a proof can be found in \cite[pp.40-42]{BoSi}: \emph{If
$\alp \in \Le[\infty] (\T)$ and $T(\alp)$ is invertible, then
\begin{itemize}
\item[1)] $\Lam(\alp)=T^{-1}(\alp)H(\alp)H(\tilde{\alp}^{-1})
=H(\alp^{-1})H(\tilde{\alp})T^{-1}(\tilde{\alp})$,
\item[2)] $T_k(\alp)B_k = P_k-P_kT(\alp)Q_k\Lam(\alp)P_k-
W_kT(\tilde{\alp}) Q_k \Lam(\tilde{\alp})W_k. $
\end{itemize}}

\medskip

\textbf{\underline{Proof} [of theorem \ref{t4}].} Throughout the
proof we assume that \linebreak $k\geq k_0$ and $1\leq j\leq k/2$.
Decompose
\begin{gather*}
 T_k^{-1}(\alp) = B_k+C_k,  \\
 C_k:=T_k^{-1}(\alp)-B_k = T_k^{-1}(\alp)(P_k-T_k(\alp)B_k).
\end{gather*}
Then
\begin{equation} \label{e5}
\begin{aligned}
 T_k^{-1}(\alp)\del_j-T^{-1}&(\alp)\del_j  =
 [T^{-1}_k(\alp)-P_kT^{-1}(\alp)P_k]\del_j+
 [P_kT^{-1}(\alp)P_k-T^{-1}(\alp)]\del_j\\
 & = [B_k+C_k-P_kT^{-1}(\alp)P_k]\del_j+
 [P_kT^{-1}(\alp)P_k-T^{-1}(\alp)]\del_j \\
 & =C_k\del_j+W_k\Lam(\tilde{\alp})W_k\del_j+
 [P_kT^{-1}(\alp)P_k-T^{-1}(\alp)]\del_j.
\end{aligned}
\end{equation}
Fix $1\leq p <q$. In order to find the parameter $c$, we will
estimate the norm of the three terms in the sum at the end of
\eqref{e5}. For this, let
\begin{gather*}
 c_1:=\sum_{n\in\Z}|n|^p|\widehat{\alp^{-1}}(n)|, \\
 c_2:=\left(\sum_{n\in\Z}|n|^{2p}|\widehat{\alp^{-1}}(n)|^2\right)^{1/2}.
\end{gather*}
Since $\alp\in W^q$ and $0\not \in \alp(\T)$, both $c_1$ and hence
$c_2$ are finite. By definition of $Q_k$ and $H(\alp^{-1})$,
\begin{align*}
 \|Q_kH(\alp^{-1})\| & \leq \|\alp^{-1}(\cdot)-\sum_{n\leq k-1}
 \widehat{\alp^{-1}}(n)\e^{in(\cdot)} \|_\infty \\
 & \leq \sum_{|n| \geq k}|\widehat{\alp^{-1}}(n)| \\
 & \leq\frac{ \sum_{|n| \geq k}|n|^p|\widehat{\alp^{-1}}(n)|}{k^p}
  \leq \frac{c_1}{k^p}.
\end{align*}
Hence by virtue of 1),
\begin{equation} \label{e5.5}
 \|Q_k\Lam(\alp)\| \leq \frac{c_1\|H(\tilde{\alp})T^{-1}(\tilde{\alp})\|}{k^p}
 \leq \frac{c_1\tilde{M}\|\alp\|_\infty}{k^p}.
\end{equation}
In a similar manner one can show
\[
 \|Q_k\Lam(\tilde{\alp})\| \leq
 \frac{c_1\tilde{M}\|\tilde{\alp}\|_\infty}{k^p}=
 \frac{c_1\tilde{M}\|\alp\|_\infty}{k^p}.
\]

We estimate the norm of $C_k$. By virtue of 2),
\begin{align*}
 C_k&=-T^{-1}_k(\alp)(T_k(\alp)B_k-P_k) \\
 & = T^{-1}_k(\alp)(P_kT(\alp)Q_k\Lam(\alp)P_k+W_kT(\tilde{\alp}) Q_k
 \Lam(\tilde{\alp})W_k).
\end{align*}
Thus
\begin{equation} \label{e6}
\begin{aligned}
 \|C_k\| &\leq \|T^{-1}_k(\alp)\|\,\|P_kT(\alp)Q_k\Lam(\alp)P_k +
 W_kT(\tilde{\alp}) Q_k\Lam(\tilde{\alp}) W_k\| \\
 &\leq M\left(\|\alp\|_\infty\|Q_k\Lam(\alp)\|+
 \|\tilde{\alp}\|_\infty\|Q_k\Lam(\tilde{\alp})\|\right) \\
 &\leq \frac{2Mc_1\tilde{M}\|\alp\|^2_\infty}{k^p}.
\end{aligned}
\end{equation}

The second term is
\begin{equation} \label{e7}
\begin{aligned}
 \|W_k\Lam(\tilde{\alp})W_k\del_j\| &\leq
 \|\Lam(\tilde{\alp})\del_{k-j+1}\| \\
 & = \|T^{-1}(\tilde{\alp})H(\tilde{\alp})H(\alp^{-1})\del_{k-j+1}\|
 \\
 & \leq \tilde{M} \|\alp\|_\infty \|H(\alp^{-1})\del_{k-j+1}\| \\
 & = \tilde{M} \|\alp\|_\infty \left(
 \sum_{n\geq k-j}|\widehat{\alp^{-1}}(n)|^2 \right)^{1/2} \\
 & \leq \tilde{M} \|\alp\|_\infty \left(
 \sum_{n\geq k/2}|\widehat{\alp^{-1}}(n)|^2 \right)^{1/2} \\
 &\leq  \frac{2^p\tilde{M} \|\alp\|_\infty \left( \sum_{n\geq
k/2}n^{2p}|\widehat{\alp^{-1}}(n)|^2\right)^{1/2}}{k^{p}} \\ &\leq
\frac{2^p\tilde{M} \|\alp\|_\infty c_2}{k^{p}}
\end{aligned}
\end{equation}
The third term is
\begin{equation} \label{e8}
\begin{aligned}
 \|[P_kT^{-1}(\alp)P_k-T^{-1}(\alp)]\del_j\| & =
 \|P_kT^{-1}(\alp)\del_j-T^{-1}(\alp)\del_j\| \\
 & = \|Q_kT^{-1}(\alp)\del_j\| \\
 & = \|Q_k[\Lam(\alp)+T(\alp^{-1})]\del_j\| \\
 & \leq \|Q_k\Lam(\alp)\del_j\|+\|Q_kT(\alp^{-1})\del_j\| \\
 & \leq \frac{c_1\|\alp\|_\infty \tilde{M}}{k^p}+
 \|Q_kT(\alp^{-1})\del_j\|\\
 & =\frac{c_1\|\alp\|_\infty \tilde{M}}{k^p} +
 \left(\sum_{n\geq k-j+1} |\widehat{\alp^{-1}}(n)|^2\right)^{1/2}
 \\ & \leq \frac{c_1\|\alp\|_\infty \tilde{M}}{k^p} +
 \left(\sum_{n\geq k/2} |\widehat{\alp^{-1}}(n)|^2\right)^{1/2}
 \\& \leq \frac{c_1\|\alp\|_\infty \tilde{M}}{k^p} +
 2^p\left(\sum_{n\geq k/2} (n/k)^{2p}|\widehat{\alp^{-1}}(n)|^2
 \right)^{1/2}
 \\& \leq \frac{c_1\|\alp\|_\infty \tilde{M}+2^pc_2}{k^p}
 \end{aligned}
\end{equation}
Hence the conclusion of theorem~\ref{t4} can be recovered from
\eqref{e5}, \eqref{e6}, \eqref{e7} and \eqref{e8} by putting
\[
 c=\tilde{M}\|\alp\|_\infty(2Mc_1\|\alp\|_\infty+2^pc_2+c_1)+2^p c_2. \find
\]

\section{The second order spectra of $H$} \label{s4}
We are now ready to state and prove the main results of this paper.
These roughly say that for all $k$ large enough,
 $\spec[2]{H;k}$ is contained in a small neighbourhood
of $\spec{H}$ with diameter of  order some negative powers of $k$.
As we previously mentioned,  applying corollary~\ref{t6}
to the decomposition \eqref{e1} will play a crucial part into the proofs.

\begin{theorem} \label{t7}
If \eqref{e2} holds, then there exists a constant $b(r)>0$
independent of $z$, such that
\[
 \spec[2]{H;k} \subset \{z\in \C\,:\, \dist{z,\mathrm{Spec}\,H} <
 b(r)/k^{2r/(2r+35)}\},
\]
for all $k$.
\end{theorem}

This automatically implies:

\begin{corollary} \label{t8}
Suppose that $r>0$ can be chosen arbitrarily large for $v$ in
\eqref{e2}. Then for all $0<q<1$ there exists a constant $b(q)>0$
independent of $z$, such that
\[
 \spec[2]{H;k} \subset \{z\in \C\,:\, \dist{z,\mathrm{Spec}\,H} <
 b(q)/k^q\},
\]
for all $k$.
\end{corollary}

The rest of this section is devoted to proving theorem~\ref{t7}.
We first introduce some further notation. Let
\[
 d(z):=\dist{z,[-2,2]}\qquad \mathrm{and} \qquad
 \tilde{d}(z):=\dist{z,\mathrm{Spec}\,H}.
\]
Then $d(z)\geq \tilde{d}(z)$. For $z\not \in [-2,2]$, we denote by
$T(\bet^{-1})$ the Toeplitz operator whose symbol is
\[
 \bet^{-1}(t;z) =(t^{-2}-2zt^{-1} +(2+z^2)-2zt+t^2 )^{-1}.
\]
Since $\bet^{-1}(\T;z)$ is an open curve and $0\not \in
\bet^{-1}(\T;z)$, then  $T(\bet^{-1})$ is invertible. Put
$\tilde{M}_2:=\|T^{-1}(\bet^{-1})\|$.

We break the proof into various steps. The aim of these steps is
to show that for all $R>0$, there exists $\tilde{b}>0$ independent
of $z$, such that $P_k(H-z)^2P_k$ is invertible provided
\[
 z\not \in \spec{H}, \qquad \tilde{d}(z)\leq R
\]
 and
\[
 k\geq \frac{\tilde{b}}{[\tilde{d}(z)]^{(2r+35)/2r}}.
\]
Since the $\spec[2]{H;k}$ are bounded uniformly in $k$ (cf.
\cite[theorem~18]{SECR}), this assertion implies theorem~\ref{t7}.
Below we assume without further mention that $z\not \in \spec{H}$
and $\tilde{d}(z)\leq R$. The various constants $b_j$ that appear
below are independent of $z$ but they may depend on $R$, $p$ or
$V$.

\underline{Step 1}: we estimate $\tilde{M}$ and $\tilde{M}_2$ for
the symbol under consideration. When $z\not\in [-2,2]$,
\[
    \zeta_{\pm}:=\frac{z\pm \sqrt{z^2-4}}{2}\not =0
\]
are such that $\zeta_+=\zeta_-^{-1}$ and $|\zeta_\pm|\not=1$. One
of these numbers is inside the unit disk whereas the other is
outside it. Denote by $\zeta$ the one that is outside and put
\[
   \beta_{\pm}(t;z):=\frac{1}{\zeta}(t^{\pm 1}-\zeta)^2.
\]
Then
\[
 \bet(t;z)=\bet_+(t;z)\bet_-(t;z),
\]
$\beta_+^{\pm 1}$ can be extended analytically inside the unit
disk and $\beta_-^{\pm 1}$ can be extended analytically outside
it. This is the so called Wiener-Hopf factorization of $\bet$, so
that (see for instance \cite[theorem~1.15]{BoSi})
\[
 T^{-1}(\bet)=T(\bet^{-1}_+)T(\bet_-^{-1}).
\]
Hence
\[
\begin{aligned}
 \tilde{M} & = \|T^{-1}(\bet)\| \\
 & \leq \|T(\bet_+^{-1})\|\,\|T(\bet_-^{-1})\|
 =\|\bet^{-1}_+\|_\infty\|\bet^{-1}_-\|_\infty \\
 & = |\zeta|^2 \sup_{t\in \T} |t-\zeta|^{-2}\sup_{t\in \T}
 |1/t-\zeta|^{-2} \\&=|\zeta|^2
 \sup|t-\zeta|^{-4} \\
 &\leq \frac{b_0}{d(z)^4}.
\end{aligned}
\]
By writing the corresponding Wiener-Hopf factorization for
$\bet^{-1}$ we obtain in a similar manner
\[
\begin{aligned}
 \tilde{M}_2 & \leq \|T(\bet_+)\|\,\|T(\bet_-)\| \\
 & = \|\bet_+\|_\infty \|\bet_-\|_\infty = \frac{1}{|\zeta|^2}
 \sup_{t\in \T} |t-\zeta|^{2}\sup_{t\in \T}
 |1/t-\zeta|^{2} \leq b_1.
\end{aligned}
\]

\underline{Step 2}: we now want to estimate the constants $c_1$
and $c_2$ in the proof of theorem~\ref{t4}, in terms of $z$.
Notice that since $\bet$ is a trigonometric polynomial, $\bet$ and
hence $\bet^{-1}$ belong to $W^p$ for all $p=1,2,\ldots$. It is
easy to show that
\[
  \frac{\ud^p}{\ud \tet^p}(\bet^{-1})(\e^{i\tet}) =
  \sum _{r=3}^{p+2} \frac{\phi_r(\tet)}{(2\cos \tet-z)^r},
\]
where the $\phi_r(\tet)$ are smooth, bounded and independent of
$z$. Then
\begin{align*}
 c_2& =
 \left(\sum_{n\in\Z}|n|^{2p}|\widehat{\bet^{-1}}(n)|^2\right)^{1/2}
  = \left(\sum_{n\not =0}|\widehat{(\bet^{-1})^{(p)}}(n)|^2\right)^{1/2}
 \\ & \leq \|(\bet^{-1})^{(p)}\|_{\Le(-\pi,\pi)}
  = b_2\left( \int_{-\pi}^{\pi} |(\bet^{-1})^{(p)}(\e^{i\tet})|^2
 \ud \tet \right)^{1/2}  \\ &\leq b_2\left( \int_{-\pi}^{\pi} \left(\sum _{r=3}^{p+2}
 \left|\frac{\phi_r(\tet)}
 {(2\cos \tet-z)^r}\right|\right)^2
 \ud \tet \right)^{1/2} \\
 & \leq b_3 \sum_{r=3}^{p+2} \frac{1}{d(z)^r} \leq
 \frac{b_4}{d(z)^{p+2}}.
\end{align*}
Similarly
\begin{align*}
 c_1 & =\sum_{n\in\Z}|n|^{p}|\widehat{\bet^{-1}}(n)| =
 \sum_{n\in\Z} \frac{|n|^{p+1}}{|n|}|\widehat{\bet^{-1}}(n)| \\
 & \leq \left(\sum_{n\not=0}\frac{1}{n^2}\right)^{1/2}
 \left(\sum_{n\not=0}|n|^{2p+2}|\widehat{\bet^{-1}}(n)|^2
 \right)^{1/2} \\
 &=  \sqrt{2\pi/6}\left(\sum_{n\not =0}|\widehat{(\bet^{-1})^{(p+1)}}(n)|^2\right)^{1/2}
  \leq \frac{b_5}{d(z)^{p+3}}.
\end{align*}

\underline{Step 3}: estimation of $M$. This only makes sense for
$k$ large enough.

\begin{lemma} \label{t10}
For all $R>0$ and $1<q<2$, there exists constants $b_6,\,b_7>0$
independent of $z$, such that $T_k(\bet)$ is invertible and
\begin{equation} \label{e20}
 \|T_k^{-1}(\bet)\| \leq \frac{b_6}{d(z)^8},
\end{equation}
provided
\[
 z\not\in [-2,2], \qquad d(z)\leq R
\]
and
\[
 k\geq \frac{b_7}{[d(z)]^q}.
\]
\end{lemma}
\Proof By virtue of \cite[lemma 2.9]{BoSi}, we know that
$T_k(\bet)$ is invertible if and only if
\[
 Q_kT^{-1}(\bet)Q_k|\mathrm{Ran}\, Q_k
\]
is invertible, and in this case
\begin{align*}
 T_k^{-1}(\bet)P_k=& \\ P_kT^{-1}(\bet)&P_k-P_kT^{-1}(\bet)Q_k
 ( Q_kT^{-1}(\bet)Q_k|\mathrm{Ran}\, Q_k)^{-1}Q_kT^{-1}(\bet)P_k.
\end{align*}

The truncation
\begin{equation} \label{e13}
  Q_kT^{-1}(\bet)Q_k|\mathrm{Ran}\, Q_k =
  Q_kT(\bet^{-1})Q_k|\mathrm{Ran}\, Q_k+Q_k\Lam(\bet)Q_k|\mathrm{Ran}\, Q_k,
\end{equation}
where
\[
 \Lam(\bet)=T^{-1}(\bet)-T(\bet^{-1})
\]
as in section \ref{s2}. The matrix of
$Q_kT(\bet^{-1})Q_k|\mathrm{Ran}\, Q_k$ is the same matrix
$T(\bet^{-1})$. Then, since $0\not \in \bet^{-1}(\T)$, the former
is invertible and
\[
 \left(Q_kT(\bet^{-1})Q_k|\mathrm{Ran}\, Q_k\right)^{-1}
\]
has the same matrix as $T^{-1}(\bet^{-1})$. Thus
\[
 \|(Q_kT(\bet^{-1})Q_k|\mathrm{Ran}\, Q_k)^{-1}\|=\|T^{-1}(\bet^{-1})\|.
\]

Let $p=1,2,\ldots$. By virtue of \eqref{e5.5},
\[
 \|Q_k\Lam(\bet)Q_k\| \leq
 \frac{c_1\tilde{M}\|\bet\|_\infty}{k^p}\, \qquad k\geq 1.
\]
According to the steps 1 and 2,
\[
 2c_1\tilde{M}\|\bet\|_\infty \tilde{M}_2 \leq
 \frac{b_8}{d(z)^{p+7}}.
\]
Then, if
\begin{equation*}
  k\geq \frac{b_8^{1/p}}{d(z)^{(p+7)/p}},
\end{equation*}
\begin{align*}
 \|(Q_kT(\bet^{-1})Q_k|\mathrm{Ran}\, Q_k)^{-1}&Q_k\Lam(\bet)Q_k\|   \\
 &\leq
 \|(Q_kT(\bet^{-1})Q_k|\mathrm{Ran}\, Q_k)^{-1}\|\,\|Q_k\Lam(\bet)Q_k\| \\
 &\leq
 \frac{\|T^{-1}(\bet^{-1})\|c_1\tilde{M}\|\bet\|_\infty}{k^p}\\
 &=\frac{c_1\tilde{M}\|\bet\|_\infty \tilde{M}_2}{k^p}\leq 1/2.
\end{align*}
By virtue of \eqref{e13}, for all such $k$,
$Q_kT^{-1}(\bet)Q_k|\mathrm{Ran}\, Q_k$ and so $T_k(\bet)$ are
invertible, and
\begin{align*}
 \|T^{-1}_k(\bet)\| & \leq
 \|T^{-1}(\bet)\|+2\|T^{-1}(\bet)\|^2\tilde{M}_2 \\
 & \leq \tilde{M}+2\tilde{M}^2\tilde{M}_2 \leq
 \frac{b_0}{d(z)^4}+\frac{2b_0b_1}{d(z)^8}.
\end{align*}
Take $p$ large enough to complete the proof. \fin

Below $M$ denotes the right hand side \eqref{e20}.

\underline{Step 4}: We are now ready to complete the proof of
theorem~\ref{t7}. For this we use corollary \ref{t6}. According to
\eqref{e1},
\[
 P_k(H-z)^2P_k =T_k(\bet)+P_kK(z)P_k,
\]
where
\[
 K(z)=H_0V+VH_0+V^2-2zV-|\del_1\rangle \langle \del_1|.
\]
It is easy to see that the matrix entries of $K(z)$ are
\begin{align*}
 \dotp{K(z)\del_n,\del_n} & = \del_1(n) +v^2(n)-2zv(n) \\
 \dotp{K(z)\del_n,\del_{n+1}} & = v(n)+v(n+1) \\
 \dotp{K(z)\del_{n+1},\del_n} & = v(n)+v(n+1) \\
 \dotp{K(z)\del_l,\del_m}& =0 \qquad \mathrm{if}\ |l-m|\geq 2.
\end{align*}
By virtue of \eqref{e2}, this ensures that $K(z)$ satisfies
conditions a)-b) of section~\ref{s3}.

Fix $p=(2r+1)/2$. According to steps 1-3, the constant $c$ of
theorem~\ref{t4} for this particular symbol is
\begin{align*}
 c&=\tilde{M}\|\bet\|_\infty(2Mc_1\|\bet\|_\infty+2^pc_2+c_1)+2^p c_2 \\
 & \leq \frac{b_{10}}{d(z)^{4}}\left(\frac{b_{11}}{d(z)^{p+11}}
 +\frac{b_{12}}{d(z)^{p+2}}
 +\frac{b_{13}}{d(z)^{p+3}}\right) + \frac{b_{14}}{d(z)^{p+2}} \\
 & \leq \frac{b_{15}}{d(z)^{p+15}} \leq
 \frac{b_{15}}{\tilde{d}(z)^{p+15}}.
\end{align*}
Also $\|K(z)\|\leq b_{16}$ and
\begin{align*}
 \|(T(\bet)+K(z))^{-1}T(\bet)\| &\leq
 \|(H-z)^{-1}\|^2\|\bet\|_\infty \\
 & \leq \frac{(2+|z|)^2}{\tilde{d}(z)^2} \leq
 \frac{b_{17}}{\tilde{d}(z)^2}.
\end{align*}
Then
\begin{align*}
  [c\sqrt{3}\|K\|+b 2^{r}\sqrt{\pi}(M+\tilde{M})]&
  \|(T(\bet)+K)^{-1}T(\bet)\| \leq \\
 & \leq \frac{b_{18}}{\tilde{d}(z)^{(2p+34)/2}}=
 \frac{b_{18}}{\tilde{d}(z)^{(2r+35)/2}}.
\end{align*}
Thus, by virtue of corollary~\ref{t6} and lemma~\ref{t10}, there
exists $\tilde{b}>0$ such that $P_k(H-z)^2P_k$ is invertible
provided
\[
 z\not\in \spec{H},  \qquad \tilde{d}(z)\leq R
\]
and
\[
  k \geq \frac{\tilde{b}}{[\tilde{d}(z)]^{(2r+35)/2r}}.
\]
This completes the proof of theorem~\ref{t7}. \fin

\medskip

In principle one can follow track of the constants
$b_0,\ldots,b_{18}$ for particular potentials. Nonetheless one
should be aware that $\tilde{b}$ increases its value faster than
$2^r$ as $r$ increases towards infinity.

\section{Numerical examples} \label{s5}

In this final section we discuss corollary~\ref{t8} for rank one
potentials from the numerical point of view. In particular we show
some computations of second order spectrum for this kind of
potentials and compare these with $\spec{H}$ which in this case
can be found explicitly. The results below extend those of
\cite[example 20]{SECR}.

Everywhere in this section, we assume that $v$
is such that
\begin{equation} \label{e15}
    v(n)= \left\{ \begin{array}{lc} a & n=j \\
    0 & \mathrm{otherwise,} \end{array} \right.
\end{equation}
for some $j=3,4,\ldots$ and  $a>2$. In the proof of the following
result we use some properties of the transfer matrix associated to
the difference equation
\begin{equation} \label{e14}
 Hu(n)=\lam u(n), \qquad \qquad u(0)=0,\, u(1)=1.
\end{equation}
See \cite{LaSi} and \cite{NaYa} for some recent results on the
relationship between spectral properties of $H$ and the transfer
matrix for slow decaying potentials.

\begin{proposition} Let $v$ be as in \eqref{e15}. Then the discrete
spectrum of $H$ consists of an isolated eigenvalue $\lam_a$, such
that
\[
 a<\lam_a<a+2/a.
\]
\end{proposition}
\Proof Since $a$ is positive, any eigenvalue $\lam_a$ of $H$
should be $\lam_a>2$. Let $\lam>2$ be a solution of the
eigenvalues equation \eqref{e14} where $u(n)$ is a sequence with
(at the moment) no constraint of growth at infinity. For all
$n\geq 1$,
\[
 \begin{vect} u(n) \\ u(n+1) \end{vect} = \begin{matr2} 0 & 1 \\
 -1 & \lam - v(n) \end{matr2} \begin{vect} u(n-1) \\ u(n)
 \end{vect}.
\]
Let
\[
 T:= \begin{matr2} 0 & 1 \\ -1 & \lam
\end{matr2}.
\]
Then for all $n\not= j$,
\begin{equation} \label{e16}
 \begin{vect} u(n) \\ u(n+1)\end{vect}= T \begin{vect} u(n-1) \\ u(n)
 \end{vect}.
\end{equation}
The Jordan decomposition for $T$ yields
\[
 T=U\begin{matr2} \mu_+ & 0 \\ 0 & \mu_- \end{matr2} U^{-1},
\]
where $\mu_\pm(\lam)\equiv \mu_\pm:=(\lam\pm\sqrt{\lam^2-4})/2$,
\[
 U=\begin{matr2} 1 & 1 \\ \mu_+ & \mu_- \end{matr2} \qquad
 \mathrm{and} \quad U^{-1}=(\mu_- -\mu_+)^{-1}
 \begin{matr2} \mu_- & -1 \\ -\mu_+ & 1 \end{matr2}.
\]
By applying inductively \eqref{e16} for $n$ from $1$ up to the
step $j-1$, we obtain
\begin{align*}
 \begin{vect} u(j-1) \\ u(j) \end{vect} & =
 T^{j-1} \begin{vect} u(0) \\ u(1) \end{vect} =
 U \begin{matr2} \mu_+^{j-1} & 0 \\ 0 & \mu_-^{j-1} \end{matr2} U^{-1}
 \begin{vect} 0 \\ 1 \end{vect} \\
 & = \frac{1}{\mu_- -\mu_+} \begin{vect} \mu_-^{j-1}-\mu_+^{j-1} \\
 \mu_-^{j}-\mu_+^{j} \end{vect}.
\end{align*}
Then
\begin{align*}
 \begin{vect} u(j) \\ u(j+1) \end{vect} & = \begin{matr2}
 0 & 1 \\ -1 & \lam-a \end{matr2} \begin{vect} u(j-1) \\ u(j) \end{vect}
 \\ & = \frac{1}{\mu_- -\mu_+} \begin{vect} \mu_-^{j}-\mu_+^{j} \\
 (\lam-a)(\mu_-^{j}-\mu_+^{j})-(\mu_-^{j-1}-\mu_+^{j-1})
 \end{vect}.
\end{align*}
Hence, for all $n>j$
\begin{equation} \label{e17}
\begin{aligned}
 \begin{vect} u(n) \\ u(n+1)\end{vect} & = T^{n-j} \begin{vect} u(j) \\ u(j+1)
 \end{vect} \\
 & = U \begin{matr2} \mu_+^{n-j} & 0 \\ 0 & \mu_-^{n-j} \end{matr2}
 U^{-1} \begin{vect} u(j) \\ u(j+1) \end{vect} \\
 & = U  \begin{vect} \mu_+^{n-j} w_1 \\ \mu_-^{n-j}w_2 \end{vect},
\end{aligned}
\end{equation}
where
\begin{gather*}
 w_1(\lam)\equiv w_1:= \frac{\mu_-u(j) - u(j+1)}{\mu_- - \mu_+} \qquad
\mathrm{and} \\
 w_2(\lam)\equiv w_2:= \frac{u(j+1)-\mu_+u(j)}{\mu_- - \mu_+}.
\end{gather*}

If $u$ is an eigenvector in $\le (\Z^+)$, then necessarily
$u(n)\to 0$ as $n\to \infty$. Hence
\[
 \mu_{+} ^{n-j}w_{1} \to 0 \qquad \mathrm{and} \qquad
 \mu_{-} ^{n-j}w_{2} \to 0
\]
as $n\to \infty$. Since $\lam >2$, then $\mu_+>1$ and $0<\mu_-<1$.
Notice that $w_1$ and $w_2$ are independent of $n$. Thus in order
for $u\in \le (\Z^+)$, necessarily $w_1=0$.

By substituting $u(j)$ and $u(j+1)$ in the identity for $w_1$,
\[
 w_1= \frac{\mu_-(\mu_-^j-\mu_+^j)-(\lam-a)
 (\mu_-^j-\mu_+^j)+(\mu_-^{j-1}-\mu_+^{j-1})}{(\mu_- - \mu_+)^2}.
\]
Notice also that $(\mu_+ / \mu_-)=\mu_+^2$. Then, a
straightforward computation shows that $w_1=0$, if and only if
$q(\mu_+)=0$ where
\[
 q(x):=x^{2j+1}-ax^{2j}-x^{2j-1}+a.
\]
It is easy to see that  $q'(x)$ vanishes, if and only if $x=0$ or
\[
 x=x_{\pm}:=\frac{2ja\pm\sqrt{(2ja)^2+4(2j+1)(2j-1)}}{4j+2}.
\]
Furthermore
\begin{itemize}
\item[1)] For all $j\geq 3$, $x_-<0$ and $x_+>1$,
\item[2)] $q(0)=a>0$,
\item[3)] $q(1)=0$ and
\item[4)] $q'(x)<0$ for $0<x<x_+$.
\end{itemize}
Hence, it is not difficult to see that $q(x)$ has only one root
$x_a$ in the interval $(1,\infty)$. Since $q(a)=a-a^{2j-1}<0$ and
\[
 q(a+1/a)=\frac{a^2+a^4+(a+1/a)^{2j}}{a+a^3}>0,
\]
necessarily $a<x_a<a+1/a$. Now,
\[
 \mu_+(\lam)=x_a,
\]
if and only if
\[
 \lam=\lam_a:=x_a+1/x_a.
\]
Therefore $\lam_a$ is the only possible eigenvalue of $H$ and
\[
  a<\lam_a<a+2/a.
\]

Finally we must show that $\lam_a$ is indeed an eigenvalue of $H$.
By virtue of \eqref{e17}, for all $n>j$
\[
 u(n)=\mu_-^{n-j}(\lam_a)w_2(\lam_a).
\]
Since $w_2(\lam_a)\not =0$ is fixed in $n$ and $0<\mu_-<1$, $u$
 is an $\le (\Z^+)$ eigenvector of $H$ as we required. \fin

\medskip

In the third column of table \ref{tb2} we show the value of
$\lam_a$ for selected $a$ and $j$. We found the roots of the
polynomial $q(x)$ using the internal algorithm ``roots'' that the
package Matlab provides. More complete data can be found without
difficulty by this method.

\begin{table}[t]
\begin{tabular} {|c|c|c|c|} \hline $a$ & $j$ & $\lam_a$
& Estimation of  $\lam_a$ using $\spec[2]{H;60}$ \\ \hline 3 & 3
&3.60362098809610& 3.60362098809609 - 0.00000008894009$i$
\\ 3 & 6& 3.60554979388607 & 3.60554979388608 -
0.00000005798745$i$\\ 3 & 9 & 3.60555127432257 &
3.60555134832088\\ 3 & 12 & 3.60555127546311 & 3.60555133471164
\\ \hline  3 & 5& 3.60553511316735 & 3.60553516591671 \\ 6 & 5 &
6.32455524824927 & 6.32455533827497 \\ 9 & 5 & 9.21954445506085 &
9.21954445506086 - 0.00000009147845$i$ \\ 12 & 5 &
12.16552506041796& 12.16552524742606
\\
 \hline
\end{tabular}
\vspace{.1in} \caption{$\lam_a$ for selected values of $a$ and
$k$.} \label{tb2}
\end{table}

\medskip

\begin{figure}[b]
\begin{picture}(250,250)(75,150) \includegraphics{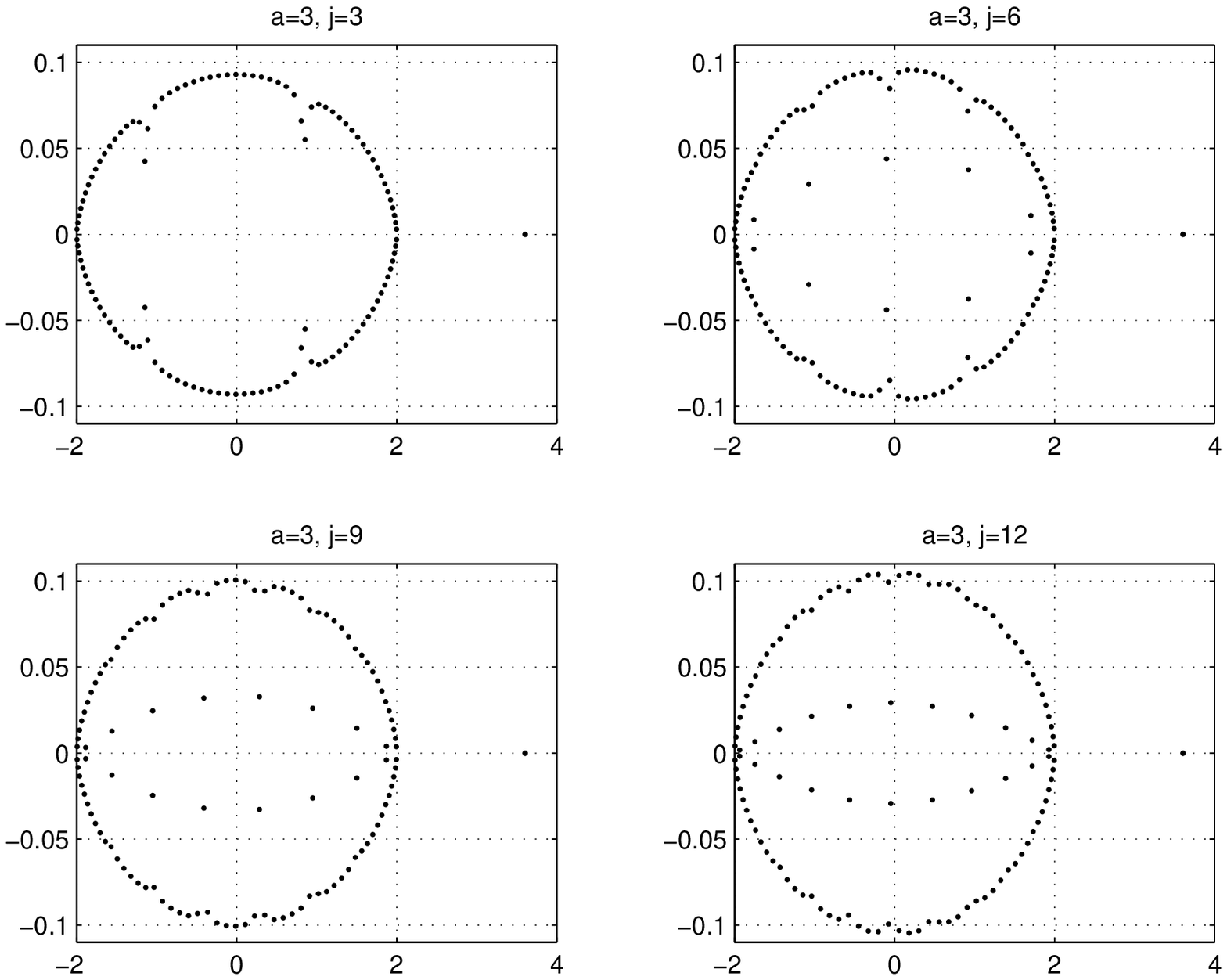}
\end{picture}
\caption{$\spec[2]{H;60}$ selected values of $a$ and $j$
corresponding to table \ref{tb2}.} \label{f4}
\end{figure}

\begin{figure}[t]
\begin{picture}(250,250)(75,150) \includegraphics{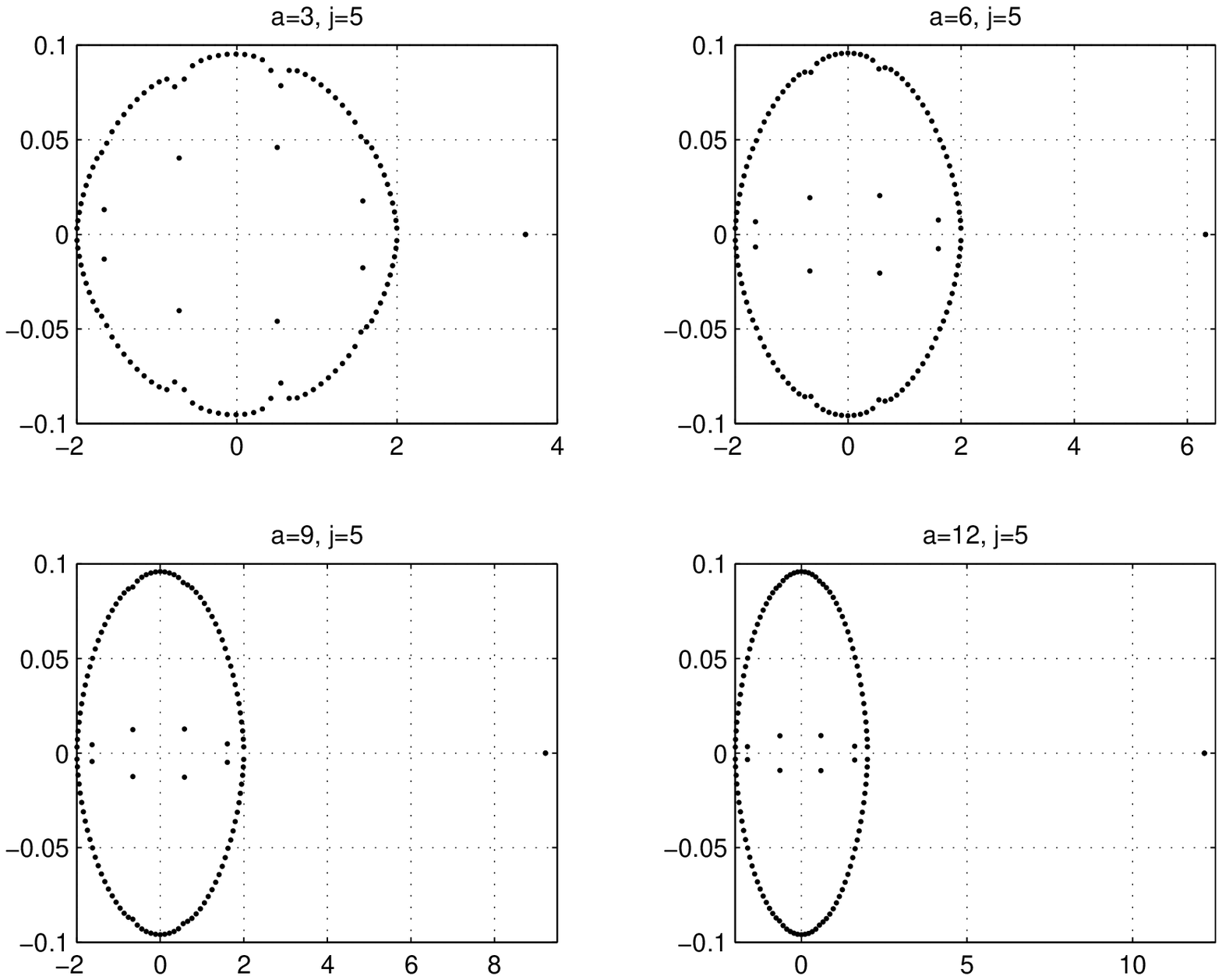}
\end{picture}
\caption{$\spec[2]{H;60}$ selected values of $a$ and $j$
corresponding to table \ref{tb2}.} \label{f5}
\end{figure}

In figures~\ref{f4} and \ref{f5} we show $\spec[2]{H;60}$ for the
first four and the last four pairs $a,j$ in table~\ref{tb2}
respectively. We found the data for this and all the other
pictures of second order spectra in this paper, by adapting in
Matlab the Maple algorithm ``reson25'' included in the appendix-a
of the electronic version of \cite{SECR}. Notice that there is
always a point in the second order spectrum which approximates the
isolated eigenvalue $\lam_a$. We reproduce the coordinates of
these points in the forth column of table~\ref{tb2}. Nonetheless
some of these values have imaginary part different from zero, in
all cases the real part of the approximation of $\lam_a$ is
accurate up to the $6^{\mathrm{th}}$ digit. This is confirmed by
further numerical experiments.

\medskip

\begin{figure}[t]
\begin{picture}(250,250)(75,150) \includegraphics{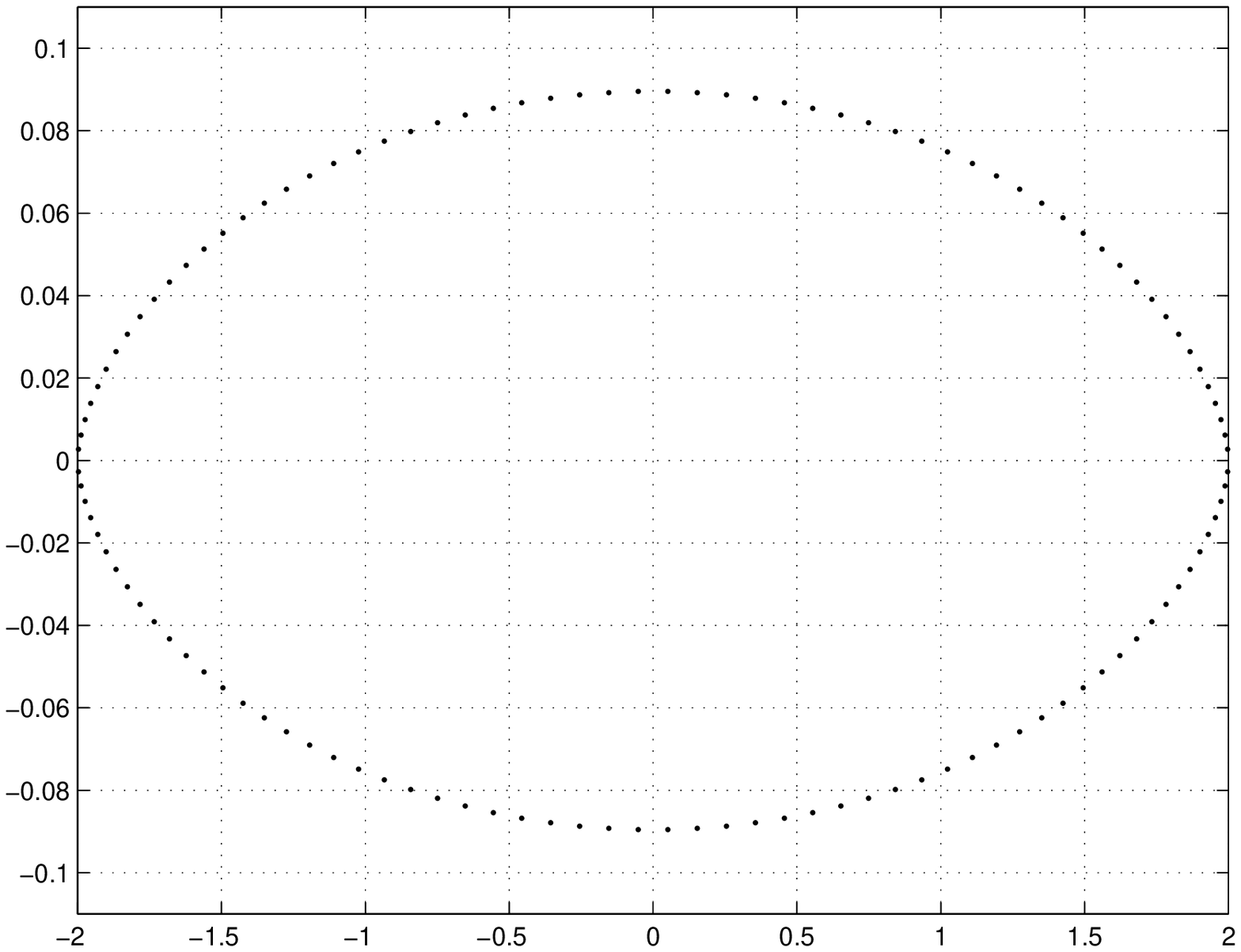}
\end{picture}
\caption{Points where $\det T_{60}(\bet)$ vanishes.} \label{f1}
\end{figure}

From figures~\ref{f4} and \ref{f5} we can also say something about
the shape of the second order spectra of $H$. Since $V$ is of rank
one, roughly speaking $\spec[2]{H;60}$ should not be too far from
the set
\[
  \{z\in \C\,:\, \det T_{60}(\bet) =0 \}
\]
which is shown in figure~\ref{f1} (the matrix is tridiagonal and
it is not very large so any of the standard computer packages can
produce the data for this picture). According to
corollary~\ref{t8} the convergence of this sets to the interval
$[-2,2]$ is $o(k^{-q})$ for all $0<q<1$ as $k\to \infty$. Compare
with figures~\ref{f4} and \ref{f5}. Notice that in the latter
there are some perturbed points very close to the real axis inside
the elliptical region which comprises most of the second order
spectrum. Further numerical experiments confirm that they increase
in number as $j$ increases. Their real part are related to the
non-real roots of $q(x)$. These non-real roots of $q(x)$ are
generalized eigenvalues which according to Davies \cite{SECR} can
be regarded as (absolute) resonances of $H$. It would be
interesting to explore further how these resonances are related to
$\spec[2]{H;k}$.

\medskip

\begin{figure}[b]
\begin{picture}(250,250)(75,150) \includegraphics{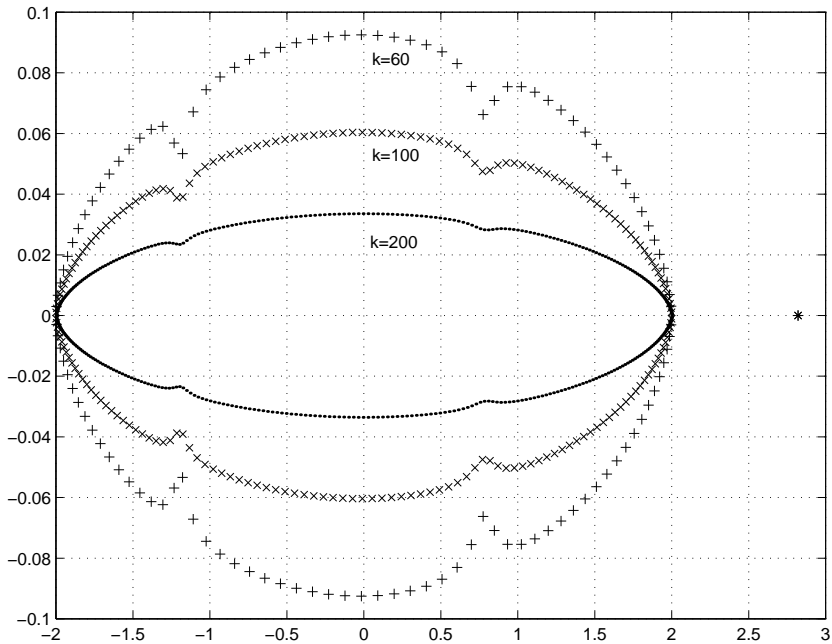}
\end{picture}
\caption{$\spec[2]{H;k}$ for $a=3$, $j=3$ and three values of
$k$.} \label{f2}
\end{figure}

\begin{figure}[t]
\begin{picture}(250,250)(75,140) \includegraphics{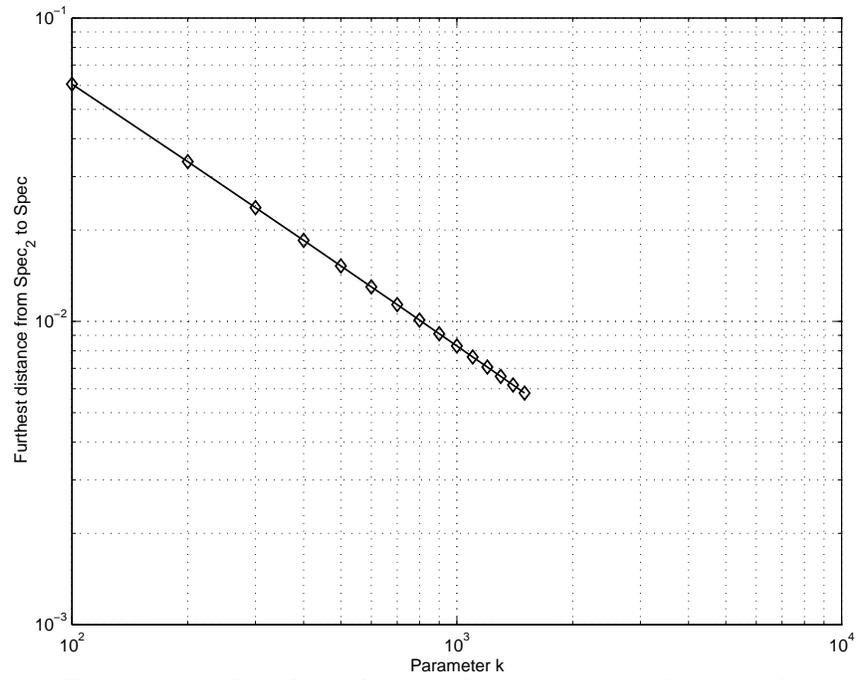}
\end{picture}
\caption{Log-log plot of the maximum distance between
$\spec[2]{H;k}$ and $\spec{H}$ for $a=3$ and $j=3$ as $k$
increases.} \label{f3}
\end{figure}

In figure~\ref{f2} we consider $\spec[2]{H;k}$ for $a=3$ and $j=3$
as $k$ increases.  Notice that the convergence to the continuous
spectrum seems to be much slower than the convergence to the
eigenvalue. The point $\phi_k\in \spec[2]{H;K}$ which is furthest
from $\spec{H}$ seems to have real part approximately equal to
$0$. Figure~\ref{f3} is a log-log graph of the imaginary part of
$\phi_k$ when $k$ varies from 100 to 1500. The slope of the line
is very close to -1. This suggests that corollary~\ref{t8} fails
for $q > 1$. In table~\ref{tb1}, the slope of the line in
figure~\ref{f3} changes from test point to test point. This
suggests that the convergence to the spectrum might be of the
order $\log^\alp(k)/k$ for some $\alp>0$.

\begin{table}
\begin{tabular} {|c|>{\PBS\centering}m{1.6in}|} \hline
k & Slope between $k$ and $k+100$.
\\ \hline100 & -0.8494 \\ 200 &  -0.8590 \\ 300 & -0.8647 \\ 400 &-0.8688
\\ 500 &-0.8719 \\600 & -0.8744 \\700 & -0.8765 \\800 &-0.8782 \\
900 & -0.8798 \\1000&  -0.8811\\1100& -0.8823\\ 1200&-0.8834  \\
1300 &-0.8844
\\1400 & -0.8853 \\ \hline
\end{tabular}
\vspace{.1in} \caption{Slope between the test points $k$ and
$k+100$ in figure~\ref{f3}.} \label{tb1}
\end{table}

\bigskip

{\samepage {\scshape Acknowledgments.} The author wishes to thank
Prof.~E.B.~Davies for his kind interest in this research and for
so many useful comments. He also wishes to thank
Prof.~S.~Marcantognini, Prof.~M.~Mor\'an, Dr.~E.~Shargorodsky and
Dr.~S.~Yakovlev for interesting discussions. Many thanks to the
referee for pointing out and suggesting how to correct a mistake
in an earlier version of this paper.}

\newpage

\end{document}